\documentclass{amsart}
\usepackage{subfigure}
\usepackage{color, amsmath, amsfonts, amssymb, epsfig, amscd, graphicx, latexsym, latexsym, amsthm, etex, synttree, tikz,comment}
\usepackage[latin5]{inputenc}
\usepackage{young}
\usepackage{mathtools}

\newtheorem{theorem}{Theorem}[section]
\newtheorem{proposition}{Proposition}[section]
\newtheorem{lemma}{Lemma}[section]

\theoremstyle{definition}
\newtheorem{definition}{Definition}[section]
\newtheorem{example}{Example}[section]
\newtheorem{remark}{Remark}[section]

\newcommand{\cellsizea}{14.142}
\newlength{\cellsza} 
\newsavebox{\cella}
\sbox{\cella}{\begin{picture}(\cellsizea,\cellsizea)
	\put(0,0){\line(1,0){\cellsizea}}
	\put(0,0){\line(0,1){\cellsizea}}
	\put(\cellsizea,0){\line(0,1){\cellsizea}}
	\put(0,\cellsizea){\line(1,0){\cellsizea}}
	\end{picture}}
\newcommand\cellifya[1]{\def\thearg{#1}\def\nothing{}%
	\ifx\thearg\nothing
	\vrule width0pt height\cellsz1 depth0pt\else
	\hbox to 0pt{\usebox{\cella} \hss}\fi%
	\vbox to \cellsza{
		\vss
		\hbox to \cellsza{\hss$#1$\hss}
		\vss}}
\newcommand\tableaux[1]{\vtop{\let\\\cr
		\baselineskip -16000pt \lineskiplimit 16000pt \lineskip 0pt
		\ialign{&\cellifya{##}\cr#1\crcr}}}

\newcommand{\cellsize}{13}
\newlength{\cellsz} 
\newsavebox{\cell}
\sbox{\cell}{\begin{picture}(\cellsize,\cellsize)
	\put(0,0){\line(1,0){\cellsize}}
	\put(0,0){\line(0,1){\cellsize}}
	\put(\cellsize,0){\line(0,1){\cellsize}}
	\put(0,\cellsize){\line(1,0){\cellsize}}
	\end{picture}}
\newcommand\cellify[1]{\def\thearg{#1}\def\nothing{}%
	\ifx\thearg\nothing
	\vrule width0pt height\cellsz depth0pt\else
	\hbox to 0pt{\usebox{\cell} \hss}\fi%
	\vbox to \cellsz{
		\vss
		\hbox to \cellsz{\hss$#1$\hss}
		\vss}}
\newcommand\tableau[1]{\vtop{\let\\\cr
		\baselineskip -16000pt \lineskiplimit 16000pt \lineskip 0pt
		\ialign{&\cellify{##}\cr#1\crcr}}}

\title[]{Lower and Upper  bounds for Nonzero Littlewood-Richardson coefficients}

\author{M\"uge Ta\c{s}k\i n}
\author{R. Bed\.i\.i G\"um\"u\c{s}}
\author{S\.inan I\c{s}\i k }
\author{M. \.ikbal Ulv\.i}
\address{ Bogaz\.I\c{c}\.I Un., Department of Mathematics,  34342 Bebek,  \.Istanbul, T\"urk\.Iye }
\address{ Bogaz\.I\c{c}\.I Un.,  Feza G\"ursey Enst.  Kand\.Ill\.I,  \.Istanbul, T\"urk\.Iye }

\thanks{\underline{This work is supported by     Bogazi\c{c}i University Research Fund Grant Number 7702 and 19082} }  
\thanks{\underline{This work is also  supported by 
T\"ubitak/1001/115F156}}


\begin{document}

	\begin{abstract} Given a skew diagram $\gamma/\lambda$, we determine a set of lower and upper bounds that a partition $\mu$ must satisfy for Littlewood-Richards coefficients $c^{\gamma}_{\lambda,\mu}>0$.   Our algorithm depends on  the characterization of  $c^{\gamma}_{\lambda,\mu}$ as the number of  Littlewood-Richardson tableau of shape $\gamma/\lambda$ and content $\mu$ and uses the (generalized) dominance order on partitions as the main ingredient.   
	\end{abstract}
	\maketitle

	\section{Introduction}
	
	\vskip.2in
	Littlewood-Richardson coefficients $c^{\gamma}_{\lambda,\mu}$ plays a vital role in many areas of mathematics. It is a well known fact that Littlewood-Richardson coefficients are the multiplicities of irreducible representations appearing  tensor products of irreducible representations of  general linear groups and symmetric group. They are  fundamental in many enumeration problems in topology and geometry since   they describe  the intersection theory of Schubert varieties of  Grassmannian.
	 
	 The calculations of  Littlewood-Richardson coefficients has been an important problem from the first time they introduced.  The first method which counts  $c^{\gamma}_{\lambda,\mu}$ as the number of Littlewood-Richardson tableaux of skew shapes is introduced in 1934 \cite{LR} but the proof is completed  thanks to the contribution of  C. Schensted \cite{Sch}, Sch\"utzenberger \cite{Schut1}\cite{Schut2} and  Knuth \cite{Knuth}. Another method  due to Berenstein and Zelevinsky \cite{BZ} calculates $c^{\gamma}_{\lambda,\mu}$ as  the number of integer points of certain polytopes. 
	 
	 The importance of these coefficients demonstrated once more with their role in the proof of Horn's  conjecture \cite{Horn}  which  basically  relates  the eigenvalues of  Hermitian matrices satisfying C = A + B through  recursively defined list of inequalities. Following  \cite{Tot} and  \cite{HR}, Klyachko  give an equivalent list  of inequalities and prove their  sufficiency in \cite{Kly} (see \cite{Fulton} for a survey).
	 
To establish a  combinatorial model encoding the list of  Klyachko,   Knutson and Tao \cite{KT1} introduce  honeycombs 	 as a new model over  the Berenstein-Zelevinsky cone. This method  brings together a lot  of important achievements \cite{KT2} \cite{KTW}: such as  the proof of saturation conjecture  and the sufficiency of the minimal set of necessary inequalities given in  \cite{Bel}.

Although the  computation of  Littlewood-Richardson coefficients  is  shown to be sharp P-complete problem by  Narayanan \cite{Nar},  the positivity of these coefficients
can be decided by a polynomial time algorithm as  first pointed out by Mulmuley and Sohoni \cite{MS} In fact, Burgisser and  Ikenmeyer \cite{BI} design such an algorithm on the base of their  characterization  of these coefficients  as the number of capacity achieving hive flows on the honeycomb graph.

It is natural to ask to what extend skew partitions can be used as a model for determining these coefficients  in terms of Klyachko's list of inequalities.  Our work presented here establishes some of the necessary  inequalities using tableaux language. Basically, for given a skew diagram  $\gamma/\lambda $ we produce a list of lower  and upper bounds   that $\mu$ has to admit if   $c^{\gamma}_{\lambda,\mu}>0$.  We note that our lower bounds (Theorem \ref{MainTheorem}) and Proposition 3.1 coincide with some of the results in \cite{McN} previously provided by P.McNamara. On the other hand since our proofs uses  different techniques, we include our approach here, to provide a more self contained exposition.  

In Section 2, we provide some background which are used in later chapters. In Section 3, we provide lower bounds (Theorem \ref{MainTheorem}) for nonzero Littlewood-Richardson coefficients whereas the upper bounds (Theorem \ref{MainTheorem2}) are given in Section 4. In the last chapter we discuss the efficiency of these upper and lower bounds.

	\section{Preliminaries} 
	\subsection{Partitions, Schur functions and Littlewood-Richardson Coefficients}
	
	For a positive integer $d$, a sequence $\lambda=(\lambda_1,\ldots \lambda_k)$ of positive integers is called a \textit{partition of $d$} if $\lambda_1\geq\ldots \geq\lambda_k>0$  and $\sum \lambda_i =d$. We write $|\lambda|=d$ in order to indicate the size of  $\lambda$.  It is natural to represent partitions with \textit{Young diagrams} which consists of left justified arrow of boxes having $\lambda_i$ boxes in the $i$-th row from the top. By \textit{a semi-standard tableau } $T$ of shape $\lambda$ we mean a labeling  the boxes of  $\lambda$ with positive integers such that labels  increasing from left to right along the rows  and strictly increasing  from top to bottom along the columns.  In this case the sequence $cont(T)=(c_1,\ldots,c_i,\ldots)$ where $c_i$ is the number of times that $i$ is used in $T$, is called \textit{the content of} $T$.
	
	For two partitions $\gamma=(\gamma_1,\gamma_2,...,\gamma_n)$  and $\lambda=(\lambda_1,\lambda_2,...,\lambda_l)$  satisfying  $\lambda\subset \gamma$, we define the \textit{ skew partition $\gamma/\lambda$} by deleting all the cells in $\gamma$ which also belongs to $\lambda$. Then \textit{a skew semi-standard tableau } of shape $\gamma/\lambda$ and its content are defined likewise. For example, 
 $$S= \tableau{&&{2}\\&{1}&{4}\\{2}&{2}}$$ 
	is a semi-standard tableau of shape  $\gamma/\lambda$ where  $\lambda=(2,1)$ and $\gamma=(3,3,2)$ and $cont(S)=(1,3,0,1)$.
	
	Any skew  semi-standard tableau $T$ of skew shape $\gamma/\lambda$ and content $cont(T)=(c_1,\ldots,c_i,\ldots)$ gives rise to the monomial $X^T=x_1^{c_1}x_2^{c_2}\ldots x_i^{c_i}$ of degree $|\gamma| -|\lambda|$. Moreover any   skew shape $\gamma/\lambda$ gives rise to  a  formal power series through the following rule
	$$s_{\gamma/\lambda}:=\sum_{T} X^T
	$$
	where $T$ runs over all skew semi-standard tableaux of shape $\gamma/\lambda$. These series are in fact among the  fundamental members of the ring of symmetric functions and they  are called   \textit{skew Schur symmetric function} indexed by  skew shapes $\gamma/\lambda$. The case  $\lambda=\varnothing$ provides  the famous  class of \textit{Schur symmetric functions}  which is  known to be a basis  for the ring of symmetric functions. Hence the product of two Schur symmetric functions can be decomposed into a sum 
	$$ s_{\lambda}s_{\mu}=\sum c_{\lambda,\mu}^{\gamma} s_{\gamma}$$
	and the coefficients $c_{\lambda,\mu}^{\gamma}$ appearing above are called  \textit{ Littlewood-Richardson (LR) coefficients}. It is also well know that any  skew Schur symmetric function $s_{\gamma/\lambda}$ can be decomposed as  
	$$s_{\gamma/\lambda}=\sum c_{\lambda,\mu}^\gamma s_{\mu}
	$$
	\subsection{ Littlewood-Richardson Rule}
	
	Littlewood-Richardson coefficients plays a vital role in many areas of mathematics. It is a well known fact that Schur symmetric functions can be realized as the characters  of the irreducible representation of general linear groups and symmetric group and  Littlewood-Richardson coefficients are the multiplicities of irreducible representations appearing  tensor products of irreducible representations of these groups. Moreover Schur symmetric functions can be realized as the fundamental classes of Schubert varieties in the chomology  of Grassmannian varieties and again Littlewood-Richardson coefficients describe  the intersection theory of Schubert varieties of these varieties  and hence fundamental in many enumeration problems in geometry. 
	
	The calculations of  Littlewood-Richardson coefficients has been an important problem from the first time they introduced.  Although the first algorithm is  given by  Littlewood and Richardson in 1934, the complete proof  is provided by Sch\"utzenberger in 1977 thanks to the contribution of  C. Schensted (1961, Insertion algorithm),  Sch\"utzenberger (1963, Jeu de taquin slides) and  Knuth (1970, Knuth Relations) along the way. Since our method in this paper highly depends  on this algorithm we provide a detailed explanation below.
	
	To start with, a word $w=w_1w_2\ldots w_i\ldots w_m$ of positive integers  is called a \textit{lattice word } if for each position $i=1,\ldots,m$, $w_i$, appears in sub-word $w_1\ldots w_i$ as many times as $w_i-1$. For example, 
	while  $w=12132$ is a  lattice word $w=1232$ is not.
	
	To continue, any skew tableau $S$ of 	shape $\gamma/\lambda$ comes with a natural word of positive integers. This word is obtained by reading the numbers in $S$ from  right to left and starting from the top row and going 
 downward direction. We call the resulting word  \textit{the reverse row  word of $S$}, denoted by $row(S)^r$. 
	
	\begin{definition} A skew tableau  $S$ of 	shape $\gamma/\lambda$ is called (LR) \textit{{Littlewood-Richardson tableau}} if $row(S)^r$ is a  lattice word. 
	\end{definition}
	As an example, consider the following tableaux of shape $\gamma/\lambda$ where  $\lambda=(2,1)$ and $\gamma=(4,3,2)$.
	
	$$S_1= \tableau{{\varnothing}&{\varnothing}&{2}&{1}\\{\varnothing}&{1}&{3}\\{1}&{2}}  \,\,\,
	S_2=\tableau{{\varnothing}&{\varnothing}&{1}&{1}\\{\varnothing}&{1}&{2}\\{3}&{3}}\,\,\,
	S_3=\tableau{{\varnothing}&{\varnothing}&{1}&{1}\\{\varnothing}&{1}&{2}\\{2}&{3}}\,\,\,
	S_4=\tableau{{\varnothing}&{\varnothing}&{1}&{1}\\{\varnothing}&{2}&{2}\\{4}&{3}}$$
	\vskip.05in
	It is easy to see that  $row(S_1)^r=123121$ and $row(S_2)^r=112133$ are not  lattice word. On the other hand $row(S_3)^r=112132$ and $row(S_4)^r=112234$ $S_4$ is not a semi-standard tableau. Hence in this example the only LR tableau is $S_3$.

	\begin{theorem} Littlewood-Richardson rule (1934-1977 by Sch\"utzenberger)
		The coefficient $c_{\lambda,\mu}^{\gamma}$ in the following equations 
		$$ s_{\lambda}s_{\mu}=\sum c_{\lambda,\mu}^{\gamma} s_{\gamma}\text{    and    } s_{\gamma/\lambda}=\sum c_{\lambda,\mu}^\gamma s_{\mu} $$
		is equal to the number of the
		{Littlewood-Richardson tableau} of  shape $\gamma/\lambda$ and  content $\mu$.
	\end{theorem}	
	
	For example $\lambda=(2,1)$ and $\gamma=(4,3,2)$ and $\mu=(3,2,1)$, the only LR tableaux that can be constructed are the followings
	$$
	\tableau{{\varnothing}&{\varnothing}&{1}&{1}\\{\varnothing}&{1}&{2}\\{2}&{3}}\,\,\,
	\tableau{{\varnothing}&{\varnothing}&{1}&{1}\\{\varnothing}&{2}&{2}\\{1}&{3}}.$$
	Hence $c_{\lambda,\mu}^{\gamma}=2$.

	\subsection{Tiling with Ballot sequences} 
	Here we provide another  characterization of- LR tableaux, that we use often later on.    To start with, observe that any skew shape $ \gamma/\lambda$ with $n$ rows and $k$ columns can be realized in $n\times k$ grid, so that any
	 cell  of $\gamma/\lambda$ lying in the $i$-th row and $j$-th column  can be identified with  coordinates $(i,j)$.  On the other hand for a tableau of shape  $ \gamma/\lambda$, a labeled cell  $c=(i,j)$ with label  $a$ can be represented by $[c,a]$.

  \begin{definition} Let $T$ be skew tableau of shape $ \gamma/\lambda$. A sequence of labeled cells   $$[(i_k,j_k),k], \ldots, [(i_2,j_2),2],[(i_1,j_1),1] $$ is called a {\it  ballot sequence } if the coordinates of these cells satisfy the following condition: 
  $$ i_r > i_{r-1} \text{ and }  j_r \leq j_{r-1} \text{ for all } 1<r\leq k. $$
That is the 
cell $(i_{r-1},j_{r-1})$ with label $r-1$ lies northeast of the cell $(i_{r},j_{r})$ with label $r$.
\end{definition}

\begin{example} One can observe that $T=\tableau{&&{1_3}&{1_2}&{1_1}\\&{1_4}&{2_2}&{2_1}\\{1_5}&{2_3}&{3_2}&{3_1}\\{4_1}}$ can be partitioned into ballot sequences listed below.  In $T$,  we use notation $i_k$ to  indicate that  the cell labeled with  $[c,i]$ belongs to  the  $k$-th ballot sequence  $B_k$.     

 $$\begin{aligned}
        & B_1=[(4,1),4],[(3,4),3], [(2,4),2],[(1,5),1] \\
        & B_2= [(3,3),3],[(2,3),2],[(1,4),1] \\
         & B_3=[(3,2),2],[(1,3),1] \\
          & B_4= [(3,1),1] \\
           & B_5=[(2,2),1] \\
    \end{aligned}
$$
\end{example}

 We remark that in the previous example the tableau $T$ that can be partitioned into ballot sequences is in fact  a LR tableau. The following lemma strengths this observation.
 
	\begin{lemma} \label{Lemma:ballot} A semi standard tableau $T$ of of shape $ \gamma/\lambda$ is an LR tableau if and only if $T$ can be partitioned in to ballot sequences.  
	\end{lemma}

	\begin{proof} We first assume that $T$ is a semi standard tableau which can be partitioned into ballot sequences. Suppose that $w_1w_2\ldots w_n$ be the reverse word of $T$ obtained by reading the labels from right to left starting from the top row. Consider the label $w_k=b$ where $1\leq k\leq n$. Now each $b$ lying in the subword $w_1w_2\ldots w_k$ lies in a different ballot sequences in $T$. Moreover in each such ballot sequence there exists cell  labeled by $b-1$ located in northeast of the cell labeled by $b$. Hence the number of  $b-1$ in the subword $w_1w_2\ldots w_k$
 is always gretaer than or equal to the number of  $b$ in the same subword. Therefore reverse row word of  $T$ is a lattice word and hence $T$ is a LR-tableau.

 We now assume that $T$ is a LR tableau and show that $T$ that can be partitioned into ballot sequences by an algorithm that construct a canonical ballot sequences. To construct  the first sequence 
  let $m$ be the largest label in $T$. Among the cells with label $m$,  let $c_m$ be the one located in the highest and the right most position. This means that $[c_m,m]$ lies in the right boundary of $\gamma/\lambda$ and that  $m$ does not label any cell in northeast   of $c_m$, except $c_m$ itself. Now the fact that the reverse word of $T$ is a lattice word yields that  at least one  label $m-1$ lies above the row of  $c_m$ and since it is the largest integer in that portion,   we can choose the one located in the right most position of the  highest possible row, say $c_{m-1}$.  In fact continuing  this way we can construct a ballot sequence $\{(c_m,m), (c_{m-1},{m-1}), \ldots,(c_2,2),(c_1,1)\}$
  in which each cell lies on the right boundary of $ \gamma/\lambda$.

 For the construction of the remaining ballot sequences, we remark that removal of the first ballot sequence from  $T$ results a tableau , say $T'$ which is still semi standard. Hence, when  we show that $T'$ is an LR tableau, the construction of ballot sequences and hence the proof of theorem is completed by induction. 
 On the other hand, this fact follows easily since the reverse word of $T'$ is obtained  by deleting the leftmost $k$ from the reverse word of $T$ for each $1\leq k\leq m$.
 
\end{proof}

	\subsection{Dominance order and its generalization} Dominance order given on the partitions, hence on Young diagrams has an important role in the areas of mathematics in which these objects are vital. The dominance order relates two partitions of the same positive integer in the following manner: Given  two partitions  $\lambda=(\lambda_1, \ldots, \lambda_k)$  and  $\mu=(\mu_1,\ldots,\mu_l) $ of  $n$, we set $\lambda \leq_{Dom} \mu$  if 
	$$  \lambda_1+\lambda_2+\ldots+\lambda_i \leq  \mu_1+\mu_2+\ldots+\mu_i  \text{ for all } 1\leq i \leq \mathrm{min}\{k,l\}.
	$$
	The method developed in this paper needs the  relation given  by  the dominance order but defined on partitions of different size. This is why we call this order the \textit{generalized dominance order.}  We note that this definition is first   introduced  in \cite{McN}.
	\begin{definition} 
		Given two partitions  $\lambda=(\lambda_1, \ldots, \lambda_k)$ of  $ n $ and $\mu=(\mu_1,\ldots,\mu_l) $ of  $m$ and $n\leq m$, we set $\lambda \leq_{GDom} \mu$ in \textit{the generalized dominance order} if 
		$$ \lambda_1+\lambda_2+\ldots+\lambda_i \leq  \mu_1+\mu_2+\ldots+\mu_i  \text{ for all } 1\leq i \leq \mathrm{min}\{k,l\}.
		$$
		
	\end{definition}
	
	\section{Lower Bounds for LR-Tableaux}
	Lower derivations  we introduce in this section plays the crucial role in determining the lower boundary for non zero LRC-coefficients.
 
		For every skew partition $\gamma/\lambda$ with $n$ rows and  $k$ columns, one can obtain a  canonical partition in the following manner:
	Let  $r_1,r_2,\ldots,r_n$ be the ordering  on the sizes  of rows in  $\gamma/\lambda$ in non-increasing fashion. Then  the partition $(r_1,r_2,\ldots,r_n)$ is called     \textit{left-top standardization} of   $\gamma/\lambda$ and denoted by 
	$$\uparrow_\leftarrow(\gamma/\lambda)=(r_1,r_2,\ldots,r_n).$$

 \begin{example} Consider $\gamma/\lambda$ where $\gamma=(5,4,4,1)$ and $\lambda=(2,1)$. Then
		$$\gamma/\lambda =\tableau{&&{}&{}&{}\\&{}&{}&{}\\{}&{}&{}&{}\\{}}
		\Longrightarrow\,\, 
		\uparrow_\leftarrow(\gamma/\lambda)=(4,3,3,1)=\tableau{{}&{}&{}&{}\\{}&{}&{}\\{}&{}&{}\\{}}.$$
	
	\end{example}

 \begin{lemma}Whenever $\gamma/\lambda\subseteq \gamma'/\lambda'$  we have
$$ \uparrow_\leftarrow(\gamma/\lambda) \subseteq \,\,  \uparrow_\leftarrow(\gamma'/\lambda') \mbox{ and hence } \uparrow_\leftarrow(\gamma/\lambda) \leq_{GDom}\,\,  \uparrow_\leftarrow(\gamma'/\lambda').$$
    
\end{lemma}

\begin{proof} Let  $r_1,r_2,\ldots,r_k$ and $r'_1,r'_2,\ldots,r'_n$ be   respectively, the orderings of the sizes of rows in  $\gamma/\lambda$ and $\gamma'/\lambda'$ in non-increasing fashion. Observe that if $r_m>r'_m$ for some $m\leq k$ then this means that   there are at least $m$ rows in $\gamma/\lambda$ with length at least $r_m$ and there are at most $m-1$ rows in $\gamma'/\lambda'$ with length at least $r_m$. But this contradicts to the fact  that  $\gamma/\lambda\subset \gamma'/\lambda'$. Hence  $r_m\leq r'_m$ for all $m\leq k$ and the result follows directly.

\end{proof}

 For a given skew partition  $\gamma/\lambda$, we denote by
 $$\gamma  /\lambda_{(i)},$$
 
 the  skew partition obtained by deleting top $i$ (possibly empty) cells in each column of $\gamma/\lambda$.
	
	Now we are ready to introduce lower derivations on $\gamma/\lambda$.

	\begin{definition} Let  $\gamma/\lambda$ be given. For $i\geq 0$, we define  \textit{$i$-th lower derivation} of $\gamma/\lambda$ by the following rule:
 $$ \partial_{(0)} (\gamma/\lambda)=\uparrow_\leftarrow(\gamma/\lambda)$$
 and for $i\geq 1$, 
		$$
		\partial_{(i)} (\gamma/\lambda)=\uparrow_\leftarrow(\gamma/\lambda_{(i)}).
		$$
	\end{definition}
	
	Observe that when $\lambda=\varnothing$, that is  $\gamma/\lambda$ is in fact a standard permutation, then  $\partial_{(1)} (\gamma/\varnothing)= \gamma=(\gamma_1,\gamma_2,\ldots)$ and 
	$\partial_{(i)} (\gamma/\varnothing)= (\gamma_i,\gamma_{i+1},\ldots)$.
	
	\begin{example} Consider $\gamma/\lambda$ given below:
		\vskip.2in
		
		$\,\,\,\,\,\,\gamma/\lambda= \tableau{&&{}&{}&{}\\&{}&{}&{}\\{}&{}&{}&{}\\{}}
		\longrightarrow 
		\partial_{(0)} (\gamma/\lambda)=\tableau{{}&{}&{}&{}\\{}&{}&{}\\{}&{}&{}\\{}}$
		\vskip.1in
		$ \gamma/\lambda_{(1)}= \tableau{&&{\varnothing}&{\varnothing}&{\varnothing}\\&{\varnothing}&{}&{}\\{\varnothing}&{}&{}&{}\\{}}
		\longrightarrow 
		\partial_{(1)}(\gamma/\lambda)=\tableau{{}&{}&{}\\{}&{}\\{}}$
		\vskip.1in
		$\gamma/\lambda_{(2)}=\tableau{& &{\varnothing}&{\varnothing}&{\varnothing}\\
			&{\varnothing}&{\varnothing}&{\varnothing}\\
			{\varnothing}&{\varnothing}&{}&{}\\{\varnothing}}
		\longrightarrow 
		\partial_{(2)}(\gamma/\lambda)=\tableau{{}&{}}$
		\vskip.1in
		$\gamma/\lambda_{(3)}=\tableau{&&{\varnothing}&{\varnothing}&{\varnothing}\\
			&{\varnothing}&{\varnothing}&{\varnothing}\\
			{\varnothing}&{\varnothing}&{\varnothing}&{\varnothing}\\{\varnothing}}
		\longrightarrow 
		\partial_{(3)}(\gamma/\lambda)=\varnothing$
		\vskip.2in
	\end{example}
	
	Now we are ready to state the following theorem.

	\begin{theorem} \label{MainTheorem} Given a skew partition $\gamma/\lambda$, if  $\mu=(\mu_1,\mu_2,\ldots,\mu_i,\ldots )$ satisfy $c_{\lambda,\mu}^{\gamma}\not=0$ then
		for all  $i\geq 1$,
		$$\partial_{(i)}(\gamma/\lambda) \leq_{GDom} (\mu_{i+1},\mu_{i+2},\ldots).  $$
		Moreover, for $\mu=\partial_{(0)}(\gamma/\lambda)$ we have  $c_{\lambda,\mu}^{\gamma}=1$, that is  
		$\partial_{(1)}(\gamma/\lambda)$ is the minimum partition satisfying $c_{\lambda,\mu}^{\gamma}\not=0$ with respect to generalized dominance order.  
	\end{theorem} 
	
	\vskip.1in

The following result is crucial in the proof of  the above theorem. To state the result we first define the following. Given a skew semi-standard tableau $T$,  let  $\#_T\colon [n] \to \mathbb{N}$ so that $\#_T(k)$ is the number of times that the label $k$ is used in $T$.

\begin{proposition} Let $T$ be a LR tableau of shape $\gamma/\lambda$ with $n$ rows and let $r_1,\ r_2,\ ... r_n$ be an ordering on the rows of 
$\gamma/\lambda$  so that their sizes satisfy $l(r_1)\geq l(r_2)\geq ...\geq l(r_n)$.  Then for all $m\leq n$, 
	$$l(r_1)+l(r_2)+ \ldots + l(r_m) \leq\, \#_T(1)+\#_T(2)+\ldots+\#_T(m)$$ 
	In particular, if $c_{\lambda,\mu}^{\gamma}\not=0$ then  $l(r_1)+l(r_2)+ \ldots+ l(r_m) \leq \,\mu_1+\mu_2+\ldots+\mu_m $ for all $m\leq n$.
\end{proposition}

\begin{proof}

By Lemma \ref{Lemma:ballot},  $T$  can be partitioned into disjoint ballot sequences. Let $P$ be the set consisting of these  disjoint  ballot sequences. Now for a fixed $m\leq n$ and for any $I\subset[m]$, let  $P^m_I$ be the set consisting ballot sequences which contains a label from   $r_i$ if and only if  $i\in I$.  Observe that such a ballot sequence can not contain a label from the row $r_j$ if $j\in [m]-I$ and may or may not contain labels from the rows  $r_{m+1},\ldots,r_n$.
It is easy to see that, for this fixed $m$
$$ P=\bigcup_{I\subset [m]} P^m_I 
$$
and this union is disjoint. Observe that if $ |I|=i$, then each ballot sequence in $P^m_I$ contains $i$ many distinct labels chosen from different rows and hence it  always contains the labels $1,2,\ldots,i$. Hence each label   $1,2,\ldots,i$  appears $ |P^m_I|$ many times in  $P^m_I$. Since $P^m_I \cap P^m_{I'}=\varnothing $ whenever $I\not = I'$ we have   
\[\#_T(1) +\#_T(2)+\ldots+ \#_T(m) \geq \sum_{i=1}^{m} \sum_{\substack{I\subset [m]\\|I|=i}} |P^m_I| + \sum_{i=2}^{m} \sum_{\substack{I\subset [m]\\|I|=i}} |P^m_I| +\ldots+ |P^m_{[m]}| =\sum_{i=1}^{m} i \sum_{\substack{I\subset [m]\\|I|=i}} |P^m_I| .   \]

On the other hand for each $1\leq j\leq m$, the row $r_j$ consists of labels belonging to different ballot sequences. Moreover each of these ballot sequence lies in  $P^m_I$ for some  $I$ containing $j$.  This shows that
$$ l(r_j)=  \sum_{\substack{I\subset [m]\\j\in I}} |P^m_I| $$
Therefore 
$$ \sum _{j=1}^m l(r_j)= \sum _{j=1}^m \sum_{\substack{I\subset [m]\\j\in I}} |P^m_I| = \sum _{j=1}^m \sum _{i=1}^m \sum_{\substack{I\subset [m]\\ |I|=i \\j\in I}} |P^m_I| = \sum _{i=1}^m  \sum_{\substack{I\subset [m]\\ |I|=i }}  \sum _{j\in I}  |P^m_I| = \sum _{i=1}^m i \sum_{\substack{I\subset [m]\\ |I|=i}} |P^m_I| $$
and hence  $\#_T(1) +\#_T(2)+\ldots+ \#_T(m) \geq \sum_{i=1}^{m} i \sum_{\substack{I\subset [m]\\|I|=i}} |P^m_I| =  \sum _{j=1}^m l(r_j)$.

\end{proof}

Now we are ready the prove the Theorem \ref{MainTheorem}.

\begin{proof} Given a skew shape $\gamma/\lambda$ with $n$ rows and $k$ columns, let $T$ be a LR-tableau of shape $\gamma/\lambda$ and with content $\mu$. Then by Lemma 4.4., for all $m\leq n$, 
$$ l(r_1)+l(r_2)+...+l(r_m) \leq \#_T(1)+\#_T(2)+....+\#_T(m)=\mu_1+\ldots+\mu_m $$ 
where   $r_1,\ r_2,\ ... ,r_n$ is the  ordering   of the rows of  $\gamma/\lambda $ so that  $l(r_1)\geq l(r_2)\geq ...\geq l(r_n)$. Now it is easy to see that  the partition  $(l(r_1), l(r_2),..., l(r_n))$ is nothing but $\partial_{(0)}(\gamma/\lambda)$.  Hence the case  $i=0$, follows directly by Lemma 4.4.

For  $i\geq 1$, first observe that removing all the cells  with label $1$,$2$,..$i-1$ from $T$ still gives a semi-standard  tableau, say $T'$ of shape $\gamma'/\lambda'$ with lowest label $i$. Taking a ballot tiling of $T$ and  deleting  $1$,$2$,..$i-1$ from each ballot sequence, we then obtained a tiling  $T'$  with ballot sequences  where the lowest label is  $i$. 
In other words $T'$ can be consider as an LR tableau of shape $\gamma'/\lambda'$  with content $(\mu_i,\ldots,\mu_n)$, where $\mu_i$ is the number of the label $i$ in $T'$. 

Hence Lemma 4.4. gives  
$$\partial_{(0)}(\gamma'/\lambda') \leq_{GDom} (\mu_i,\mu_{i+1},\ldots,\mu_n) $$
Now the theorem follows when we show that 
$\partial_{(i)}(\gamma/\lambda) \leq_{GDom} \partial_{(0)}(\gamma'/\lambda')$.

Recall  that 
$\partial_{(i)}(\gamma/\lambda)= \uparrow_\leftarrow(\gamma/\lambda_{(i-1)})
$
where $\gamma/\lambda_{(i-1)}$ is obtained by deleting top $i-1$ cells in each column of $\gamma/\lambda$ whenever possible. 

On the other hand, 
$\partial_{(0)}(\gamma'/\lambda')=\uparrow_\leftarrow(\gamma'/\lambda')
$
where $\gamma'/\lambda'$ is the shape of the tableau $T'$ obtained from $T$ by removing the cells  labeled by  $1$,$2$,..,$i-1$. Observe that these cells removed can be located only in  the first $i-1$ rows of each column of  $\gamma/\lambda$. That is 
$$\gamma/\lambda_{(i-1)}\subset \gamma'/\lambda'
$$
and by Lemma 3.1  this proves $\partial_{(i)}(\gamma/\lambda)= \uparrow_\leftarrow(\gamma/\lambda_{(i-1)}) \leq_{GDom}\,\,  \uparrow_\leftarrow(\gamma'/\lambda')=\partial_{(0)}(\gamma'/\lambda')$ as required.

\vskip.1in
Now we will show that $c_{\lambda,\mu}^{\gamma}=1$ for  $\mu=\partial_{(0)}(\gamma/\lambda)=(l(r_1), l(r_2),..., l(r_n))$. Observe that  $\gamma/\lambda$ has $n$ rows, each of which  has at least $l(r_n)$ cells. Hence for each $1\leq i\leq n$, labeling right most $l(r_n)$ cells of $i$-th row by $i$, we construct  $l(r_n)$ many ballot sequences in $\gamma/\lambda$ consisting of the labels $1,2,\ldots,n$.

Now, assuming that $k$ is the largest index such that $l(r_k)-l(r_n)>0$, we then have  the reduced content $\mu'=(l(r_1)-l(r_n),\ldots,l(r_k)-l(r_n))$. On the other hand  removing these $l(r_n)$ ballot sequences from $\gamma/\lambda$ we  have left with a sub skew shape $\gamma'/\lambda'$.   Here  
$\gamma'/\lambda'$ has also $k$ rows, each of which has  at least $l(r_k)-l(r_n)$ cells and by repeating the same process above  we can construct $l(r_k)-l(r_n)$ many ballot sequences consisting of the labels $1,2,\ldots,k$. Clearly,  continuation of  this process terminates  $\gamma/\lambda$ and the content $\mu $ completely.  Therefore $\gamma/\lambda$ can be tiled with ballot sequences by using the labels coming from with content $\mu=\partial_{(0)}(\gamma/\lambda)$. Denoting the resulting tableau with $T$ we then deduce that   $c_{\lambda,\mu}^{\gamma}\geq 1$.

Now we will show that $c_{\lambda,\mu}^{\gamma}\not > 1$. First observe that if $S$ is an LR tableau of $\gamma/\lambda$ and content $\mu=(l(r_1), l(r_2),..., l(r_n))$ then $S$ must be partitioned into ballot sequences among which $l(r_n)$ of them consisting of  $1,2,\ldots, n$. Since  $\mu$ has $n$ rows then each label $i$ must label $l(r_n)$ cells in the $i$-th row. But using this reasoning, it is easy to see that the tableau $S$ must be tiled in the same way as the tableau $T$ above. Therefore $S=T$ and hence $c_{\lambda,\mu}^{\gamma}= 1$.
\end{proof}

\section{Upper Bounds for LR-tableaux}

Similar to the lower derivations, upper derivations that we introduce in this section plays the crucial role in determining the upper boundary for non zero LRC-coefficients.

For every skew partition $\gamma/\lambda$ with $n$ rows and  $k$ columns, one can obtain a  canonical partition in the following manner:
Let  $c_1,c_2,\ldots,c_k$ be the ordering  of the sizes of columns in  $\gamma/\lambda$  in non-increasing fashion. Then the transpose of the partition $(c_1,c_2,\ldots,c_k)$ is called \textit{top-left standardization} of   $\gamma/\lambda$ and denoted by 
$$^\leftarrow\uparrow(\gamma/\lambda)=(c_1,c_2,\ldots,c_k)^t.$$

\begin{lemma} Whenever $\gamma/\lambda\subseteq\gamma'/\lambda'$  we have
$$ ^\leftarrow\uparrow(\gamma/\lambda) \subseteq\,\,  ^\leftarrow\uparrow(\gamma'/\lambda') \mbox{ and hence }  ^\leftarrow\uparrow(\gamma/\lambda) \leq_{GDom}\,\,  ^\leftarrow\uparrow(\gamma'/\lambda').$$
    \end{lemma}

\begin{proof}
It is easy to see that  $(\uparrow_\leftarrow((\gamma/\lambda)^t))^t=^\leftarrow\uparrow(\gamma/\lambda)$. Now $\gamma/\lambda\subseteq\gamma'/\lambda'$ implies that $(\gamma/\lambda)^t\subseteq (\gamma'/\lambda')^t$. Hence  by Lemma 3.1, we have 
$$\uparrow_\leftarrow((\gamma/\lambda)^t) \subseteq\,\,  \uparrow_\leftarrow((\gamma'/\lambda')^t) \mbox{ and hence } (\uparrow_\leftarrow((\gamma/\lambda)^t))^t\subseteq\,\,(\uparrow_\leftarrow((\gamma'/\lambda')^t))^t.$$ Therefore the result follows directly.
\end{proof}

\begin{example} Consider $\gamma/\lambda$ where $\gamma=(5,4,4,1)$ and $\lambda=(2,1)$. Then
	$$\gamma/\lambda =\tableau{&&{}&{}&{}\\&{}&{}&{}\\{}&{}&{}&{}\\{}}
	\,\, \text{  and  }\,\, 
	^\leftarrow\uparrow(\gamma/\lambda)=(5,4,2)=\tableau{{}&{}&{}&{}&{}\\{}&{}&{}&{}\\{}&{}}$$
\end{example}

The following result is crucial in the construction of upper bounds for LR-tableaux. 

		\begin{proposition}
		    
	Given a skew partition $\gamma/\lambda$, if  $\mu=(\mu_1,\mu_2,\ldots,\mu_j)$ satisfies $c_{\lambda,\mu}^{\gamma}\not=0$ then
 $$\mu\leq_{GDom} \,^\leftarrow\uparrow(\gamma/\lambda). $$
Moreover, for $\nu= ^\leftarrow\uparrow(\gamma/\lambda)=(\nu_1,\nu_2,\ldots)$, we have $c_{\lambda,\nu}^{\gamma}=1 $.
\end{proposition}
	\begin{proof} Let $\nu= ^\leftarrow\uparrow(\gamma/\lambda)=(\nu_1,\nu_2,\ldots)$. Observe that $\nu_i$ is the number of columns in $\gamma/\lambda$ with size $\geq i$. Let $T$ be an LR-tableau with content $\mu$. Since the label $1$ may appear only in the first cell (from top) of a column then $\mu_1$ is less than or equal to the number of columns with size $\geq 1$. That is $\mu_1\leq \nu_1$. Similarly  we have $\mu_1+\mu_2\leq \nu_1+\nu_2$, since label $1$ and label $2$ may appear only in the first or second cell of a column and since  $\nu_1+\nu_2$ is the number of cells lying in the first or second place in a column of $\gamma/\lambda$. Continuing this way we can deduce that $\mu_1+\mu_2+\ldots+\mu_i\leq\nu_1+\nu_2+\ldots+\nu_i$,
	for each $i$. That is  $\mu\leq_{GDom} \, ^\leftarrow\uparrow(\gamma/\lambda)$.
	
	Tho show  show $c_{\lambda,\nu}^{\gamma}=1$, first observe that  labelling by $i$, each cell lying in the $i$-th place (from top) of a column in  $\gamma/\lambda$, we obtained  a semi-standard tableau, say $T$, of shape $\gamma/\lambda$ and content $\nu$. In fact,  this is a $LR$-tableau since in each column of size $k$ gives us a ballot sequence. Hence, $c_{\lambda,\nu}^{\gamma}\not=0$. 
	
	Now, we will show the labelling that is described above is unique with given $\gamma/\lambda$ and $\nu$. Assume that $S$ is a LR-tableau with shape $\gamma/\lambda$ and content $\nu$. Since $\nu_1$ is the number of cells which lies in the first position in a column and since $\nu_1$ is the number of label $1$ in $S$, we see that  those cells lying in the first position of each column must be labeled by $1$. By induction we may assume that, in $S$, the top $i$ positions  in each column of size $\geq i$  must be labeled  by $1,2,\ldots,i$. Now since $\nu_{i+1}$ is the number of cells which lies in the $(i+1)$-th position in a column of size $\geq i+1$ and since $\nu_{i+1}$ is  the number of label $i+1$ in $S$, we see that  those cells lying in the $(i+1)$-th position of each column of size $\geq i+1$ should be labeled by $i+1$. Hence this argument shows that $S$  is nothing but the tableau $T$ constructed in the previous paragraph.
	\end{proof}

Let  $\gamma/\lambda$ be a skew partition with $n$ rows and  $k$ columns. For any  sequence $1\leq i_1<\ldots<i_j\leq n$, we define $$\gamma  /\lambda^{(i_1,\ldots,i_j)}$$
to be    the  skew partition obtained by deleting the rows of $\gamma/\lambda$ indexed by $ i_1,\ldots,i_j$ and concatenating remaining rows vertically. Now we are ready to define upper derivations on skew partitions.

\begin{definition} Let  $\gamma/\lambda$ be a skew partition with $n$ rows and  $k$ columns.  
For $j\geq 0$, we define $j$-th  upper derivation of $\gamma/\lambda$ by the following rule: if $j=0$,
$$ \partial^{(0)} (\gamma/\lambda):= \,^\leftarrow\uparrow(\gamma/\lambda) $$
and for $j\geq 1$,  any sequence $1\leq i_1<\ldots<i_{j-1}\leq n$ defines a $j$-th derivation of  $\gamma/\lambda$ by the following rule 
	$$
	\partial^{(i_1,\ldots,i_j)} (\gamma/\lambda):=\,^\leftarrow\uparrow(\gamma/\lambda^{(i_1,\ldots,i_{j})}).
	$$
\end{definition}

\begin{example} Consider $\gamma/\lambda$ given below:
	\vskip.2in
	
	$\,\,\,\,\,\,\gamma/\lambda= \tableau{&&{}&{}&{}\\&{}&{}&{}\\{}&{}&{}&{}\\{}}
	\longrightarrow 
	\partial^{(0)} (\gamma/\lambda)=\tableau{{}&{}&{}&{}&{}\\{}&{}&{}&{}\\{}&{}}$
	\vskip.1in
	$  \tableau{&&{}&{}&{}\\&{}&{}&{}\\{\varnothing}&{\varnothing}&{\varnothing}&{\varnothing}\\{}}
	\longrightarrow 
	\gamma/\lambda^{(3)}=\tableau{&&{}&{}&{}\\&{}&{}&{}\\{}}
	\longrightarrow  \partial^{(3)}(\gamma/\lambda)=\tableau{{}&{}&{}&{}&{}\\{}&{}}$
	\vskip.1in
	$\tableau{&&{\varnothing}&{\varnothing}&{\varnothing}\\&{}&{}&{}\\{\varnothing}&{\varnothing}&{\varnothing}&{\varnothing}\\{}}
	\longrightarrow 
	\gamma/\lambda^{(1,3)}=\tableau{&{}&{}&{}\\{}}
	\longrightarrow 
	\partial^{(1,3)}(\gamma/\lambda)=\tableau{{}&{}&{}&{}}$
	\vskip.1in
	$\tableau{&&{\varnothing}&{\varnothing}&{\varnothing}\\
		&{\varnothing}&{\varnothing}&{\varnothing}\\
		{\varnothing}&{\varnothing}&{\varnothing}&{\varnothing}\\{}}
	\longrightarrow 
	\gamma/\lambda^{(1,2,3)}=\tableau{
	{}}
	\longrightarrow 
	\partial^{(1,2,3)}(\gamma/\lambda)=\tableau{{}}$
	\vskip.1in
	and $
	\partial^{(1,2,3,4)}(\gamma/\lambda)=\varnothing$
	\vskip.2in
	
\end{example}

Now we are ready to state the  theorem on upper bounds of $\gamma/\lambda$.

\begin{theorem} \label{MainTheorem2} Given a skew partition $\gamma/\lambda$ with $n$ rows and $k$ columns, if  $\mu=(\mu_1,\mu_2,\ldots,\mu_i,\ldots )$ satisfy $c_{\lambda,\mu}^{\gamma}\not=0$ then 
$$ (\mu_{1},\mu_{2},\ldots) \leq_{GDom}    \partial^{(0)}(\gamma/\lambda)   $$
and for all $1\leq j \leq n$ and for all sequence  $1\leq i_1<\ldots<i_j\leq n$ we have 
	$$ (\mu_{j+1},\mu_{j+2},\ldots) \leq_{GDom}    \partial^{( i_1,\ldots,i_j)}(\gamma/\lambda)   $$
	Moreover, for $\mu=\partial^{(0)}(\gamma/\lambda)$ we have  $c_{\lambda,\mu}^{\gamma}=1$, that is  
	$\partial_{(0)}(\gamma/\lambda)$ is the maximum partition satisfying $c_{\lambda,\mu}^{\gamma}\not=0$ with respect to generalized dominance order.  
\end{theorem} 
\subsection{Row Bumping}  We first explain row bumping algorithm  that we use as a main tool in proving the above theorem.  Let $T$ be an LR tableau of shape $\gamma/ \lambda$ and $B_1,\ldots, B_n$ be ballots sequences which tile $T$. Here any label $[c,j]$ in $B_k$ is represented with sub-index $j_k$ to indicate that $[c,j]$ belongs to $B_k$. 

Basically, row bumping algorithm    bumps up all labels  lying in the $i$-th row, starting from the right most label say $ j_k$ and proceeding to the left. On the other hand, to bump up the label $ j_k$ we traces its ballot  sequence $B_k$ in $T$ and  check  whether the replacement of   $ (j-1)_k$ with $ j_k$ still gives a semi-standard tableau. If it does then $ (j-1)_k$ is replaced with $ j_k$ and bumping up process now continues with $ (j-1)_k$ along $B_k$ again.  

If the replacement of   $ (j-1)_k$ with $ j_k$  does not give  a semi-standard tableau then we proceed as follows: 
If the replacement of $ (j-1)_k$ with $ j_k$  disturb   the  strict increase condition in  columns then this means that in the original tableau $T$, $ (j-1)_k$ sits above $j_r$ where $j_r$ lies in the ballot sequence $B_r$. In this case we interchange the labeled cells $ (j-1)_k,\ldots,2_k,1_k$ of $B_k$  with $(j-1)_r,\ldots,2_r,1_r$ of $B_r$ to get another tiling of $T$ and we continue to bumping up $ j_k$ along new $B_k$ in this new tiling.

On the other hand,  if the replacement of $ (j-1)_k$ with $ j_k$ does not disturb  the  strict increase condition in  columns  but  disturbs the non-decrease condition in rows, then this means that to the right of $(j-1)_k$ there is at least one cell labeled by $(j-1)$. Among these, labeled cells let $(j-1)_s$ be the one in the right most position. Now  $(j-1)_s$ lies in the ballot sequence $B_s$ and  we interchange   $ (j-1)_k,\ldots,2_k,1_k$ of $B_k$  with $(j-1)_s,\ldots,2_s,1_s$ of $B_s$ to get another tiling of $T$. Then we continue to bumping up $ j_k$ along new $B_k$ in this new tiling. Note that in this case, the replacement of $ (j-1)_k$ with $ j_k$  disturb neither  the  strict increase condition in  columns nor the non-decrease condition in rows.

Hence, applying the process above whenever needed we arrive in a new tiling of $T$ in which, along the new ballot sequence $B_k$,  $ j_k$ replaces $ (j-1)_k$,  $ (j-1)_k$ replaces $ (j-2)_k$, so on so forth and bumps  $1_k$ out from the tableau $T$. We denote this new tableau by
$$T^{\uparrow j_k}.
$$
 Note that, 
 \begin{enumerate}
     \item all ballot sequence in last tiling of $T$, except $B_k$, remains same in  $T^{\uparrow j_k}$,
     \item new $B_k$ in $T^{\uparrow j_k}$ can be considered as a ballot sequence in which the smallest label is $2_k$ and this label sits in the top cell of a column,
     \item for each  $ 2_k\leq r_k \leq j_k$ in new $B_k$, if there is a label sits above $r_k$, this label must be $\leq r_k-2$ and if there is a label sits to the left of $r_k$, this label must be $\leq r_k-   1$. 
 \end{enumerate}
 
   Now to achieve row bumping on the   $i$-th row of $T$, we start by bumping the rightmost labeled cell say $j_k$ in $i$-th row to find $T^{\uparrow j_k}$ and   continue with  the next left labeled cell, say $j'_s$, to find $(T^{\uparrow j_k})^{\uparrow j'_s}$,   until all the cells in this row are exhausted.  Here observe that because of the property (3) above, the ballot sequence $B_k$ of $T^{\uparrow j_k}$ remains same in $(T^{\uparrow j_k})^{\uparrow j'_s}$ that is the cell removed has also label $1$. 

When we bump all the cells lying in the $i$-th row of $T$, the tableau we get, in general, has not a skew shape but,   concatenating $(i+1)$-th row to the  $(i-1)$-row vertically gives  a semi-standard tableau, denoted by 
$$ T^{(i)} \text{ of shape  } \gamma/\lambda^{(i)} \text{ and content } (\mu_1-r_i,\mu_2,\mu_3,\ldots,).$$ 
where  $r_i$ is the size of the $i$-th row.

 \begin{example}
For an example,  with consider the first tableau, $T$, below with shape $\gamma/\lambda$ with  content $\mu=(7,7,5,4,1)$. It is easy to observe that $T$ can be tiled with seven ballot sequences $B_1$ through $B_7$ where the number $j$ in $B_i$ is indicated with sub-index $j_i$. Hence $T$ is an LR-tableau and that $(7,7,5,4,1) \leq_{GDom} \, ^\leftarrow\uparrow(\gamma/\lambda)=(10,8,5,1)$. 
  \vskip .1in $$ \tableau{&&&&&&{1_4}&{1_2}&{1_1}&{1_6}&{1_7}\\&&&&&{\bf 1_3}&{2_2}&{2_1}&{2_6}&{2_7}\\&&&&{1_5}&{\bf2_3}&{3_1}\\&&&&{\bf 3_3}&{3_2}\\
 &&{2_5}&{2_4}&{4_2}&{4_1}\\&&{3_4}&{\bf 4_3}\\&{3_5}&{4_4}&{5_3} }  \tableau{&&&&&&{1_4}&{\bf1_3}&{1_1}&{1_6}&{1_7}\\&&&&&{1_2}&{\bf2_3}&{2_1}&{2_6}&{2_7}\\&&&&{1_5}&{2_2}&{3_1}\\&&&&{ 3_2}&{\bf3_3}\\
 &&{2_5}&{2_4}&{4_2}&{4_1}\\&&{3_4}&{\bf 4_3}\\&{3_5}&{4_4}&{5_3} }   \tableau{&&&&&&{1_4}&{1_1}&{\bf1_3}&{1_6}&{1_7}\\&&&&&{1_2}&{2_1}&{\bf2_3}&{2_6}&{2_7}\\&&&&{1_5}&{2_2}&{\bf3_3}\\&&&&{ 3_2}&{3_1}\\
 &&{2_5}&{2_4}&{4_2}&{4_1}\\&&{3_4}&{\bf 4_3}\\&{3_5}&{4_4}&{5_3} } 
 \,\,$$  
   \vskip .1in
We now explain the process on the labeled cell $4_3$ which belong to $B_3$. The algorithm aims  to  bump up $4_3$ in its ballot sequence, that is $4_3$ takes place of $3_3$, as long as this replacement still gives strictly increasing numbers in the columns. In  our example,  $3_3$ sits above the label  $4_2$, and hence    $4_3$ can not take the place of $3_3$. To overcome this problem we interchange  the places of $3_3,2_3,1_3$ in $B_3$ with $3_2,2_2,1_2$ in $B_2$ to get another ballot tiling of $T$, as it is shown in the second tableau above. On the other hand in the second tiling of $T$,  $3_3$ sits above $4_1$, so again $4_3$ can not take place of  $3_3$ and we continue with interchanging the places of  $3_3,2_3,1_3$ in $B_3$ with $3_1,2_1,1_1$ in $B_1$, to 
  arrive a new tiling  of $T$, shown in the third tableau above.

  At this point, we can  put $4_3$ in the place of  $3_3$, so we check now whether  $3_3$ can take place of $2_3$ in its ballot sequence.  Now we see that placing  $3_3$ in the place of $2_3$ does not disturb   the  strict increase condition in  columns but it disturbs the non-decrease condition in rows. To solve this problem we interchange the  places of  $2_3,1_3$  of $B_3$ with  $2_7,1_7$ of $B_7$ as it is shown the second tableau below. In this new tiling, $3_3$ takes place of $2_3$ and also then $2_3$ takes place of $1_3$. as  shown in the first tableau below.  Then by doing these replacements, we obtain  the second tableau below which is  $T^{\uparrow 4_3}$. 
  
  \vskip .1in
 $$ \tableau{&&&&&&{1_4}&{1_1}&{1_7}&{1_6}&{\bf 1_3}\\&&&&&{1_2}&{2_1}&{2_7}&{2_6}&{\bf 2_3}\\&&&&{1_5}&{2_2}&{\bf 3_3}\\&&&&{ 3_2}&{3_1}\\
 &&{2_5}&{2_4}&{4_2}&{4_1}\\&&{3_4}&{\bf 4_3}\\&{3_5}&{4_4}&{5_3} } \, \, \, \, \, \, \, \, \, \, \, \, \, \, \, \,
T^{\uparrow 4_3}= \tableau{&&&&&&{1_4}&{1_1}&{1_7}&{1_6}&{\bf 2_3}\\&&&&&{1_2}&{2_1}&{2_6}&{2_7}&{\bf 3_3}\\&&&&{1_5}&{2_2}&{\bf 4_3}\\&&&&{ 3_2}&{3_1}\\
 &&{2_5}&{2_4}&{4_2}&{4_1}\\&&{3_4}&\\&{3_5}&{4_4}&{5_3} }
 $$ 
  \vskip .1in
Now one can check easily that continuing bumping up process on the second labeled cell $3_4$ in the sixth row, gives the following sequence of tiling on the tableau $T^{\uparrow 4_3} $.
  \vskip.1in $$
  \tableau{&&&&&&{\bf1_4}&{1_1}&{1_7}&{1_6}&{2_3}\\&&&&&{1_2}&{2_1}&{2_6}&{2_7}&{ 3_3}\\&&&&{1_5}&{2_2}&{ 4_3}\\&&&&{ 3_2}&{3_1}\\
  	&&{2_5}&{\bf2_4}&{4_2}&{4_1}\\&&{\bf3_4}&\\&{3_5}&{4_4}&{5_3} }
  \tableau{&&&&&&{1_1}&{\bf1_4}&{1_7}&{1_6}&{2_3}\\&&&&&{1_2}&{2_1}&{2_6}&{2_7}&{ 3_3}\\&&&&{1_5}&{2_2}&{ 4_3}\\&&&&{ 3_2}&{3_1}\\
  	&&{2_5}&{\bf2_4}&{4_2}&{4_1}\\&&{\bf3_4}&\\&{3_5}&{4_4}&{5_3} }
 \tableau{&&&&&&{1_1}&{1_6}&{1_7}&{\bf 1_4}&{2_3}\\&&&&&{1_2}&{2_1}&{2_6}&{2_7}&{ 3_3}\\&&&&{1_5}&{2_2}&{ 4_3}\\&&&&{ 3_2}&{3_1}\\
 	&&{2_5}&{\bf2_4}&{4_2}&{4_1}\\&&{\bf3_4}&\\&{3_5}&{4_4}&{5_3} }
 $$ 
   \vskip .1in
  In the last tiling above, $3_4$ bumps up $2_4$ and then bumps up $1_4$ without disturbing the  strict increase condition in  columns and  the non-decrease condition in rows, so we get $(T^{\uparrow 4_3})^{\uparrow 3_4}$ as shown below.  Now  concatenating $7$- the row to the $5$-th row of $(T^{\uparrow 4_3})^{\uparrow 3_4}$ vertically, we obtained the semi-standard  tableau, denoted by  $T^{(6)}$, with shape $ \gamma/\lambda^{(6)}$.   \vskip .1in
   $$ (T^{\uparrow 4_3})^{\uparrow 3_4} = \tableau{&&&&&&{1_1}&{1_6}&{1_7}&{2_4}&{2_3}\\&&&&&{1_2}&{2_1}&{2_6}&{2_7}&{3_3}\\&&&&{1_5}&{2_2}&{4_3}\\&&&&{ 3_2}&{3_1}\\
   	&&{2_5}&{3_4}&{4_2}&{4_1}\\&\\&{3_5}&{4_4}&{5_3} } \,\, \,\, \,\,\,\,\,\,\  T^{(6)}=\tableau{&&&&&&{1_1}&{1_6}&{1_7}&{2_4}&{2_3}\\&&&&&{1_2}&{2_1}&{2_6}&{2_7}&{3_3}\\&&&&{1_5}&{2_2}&{4_3}\\&&&&{ 3_2}&{3_1}\\
   	&&{2_5}&{3_4}&{4_2}&{4_1}\\ &{3_5}&{4_4}&{5_3} } $$  \vskip .1in
\end{example}   

\begin{remark} The proof of Theorem \ref{MainTheorem2} relies on the following observation. In the above example,  if we  delete all the cells labeled with $1$ in $T^{(6)}$ we obtained a tableau say $T'$ with shape,  say $\gamma'/\lambda'$, and  content $(0,7,5,4,1)$.
  \vskip .1in $$ T'=\tableau{&&&&&&&&&{2_4}&{2_3}\\&&&&&&{2_1}&{2_6}&{2_7}&{3_3}\\&&&&&{2_2}&{4_3}\\&&&&{ 3_2}&{3_1}\\
 &&{2_5}&{3_4}&{4_2}&{4_1}\\&{3_5}&{4_4}&{5_3} }
 $$ 
   \vskip .1in
  It is easy to see that, by restricting the labels from $\{1,2,3,\ldots, \}$ to $\{2,3,\ldots, \}$,  $T'$ can be considered as LR-tableau with content $(7,5,4,1)=(\mu_2,\mu_3,\mu_4,\mu_5)$. Hence, by Proposition 3.1 we have 
 $$(\mu_2,\mu_3,\mu_4,\mu_5) \leq_{GDom} \, ^\leftarrow\uparrow(\gamma'/\lambda'). 
 $$
 \vskip.1in
 Now, since $\gamma'/\lambda' \subset  \gamma/\lambda^{( 6)}$, by Lemma 4.1, $^\leftarrow\uparrow(\gamma'/\lambda') \leq_{GDom} \, ^\leftarrow\uparrow(\gamma/\lambda^{( 6)}) =\partial^{( 6)}(\gamma/\lambda), $
$$\text{ and hence  } 
  (\mu_2,\mu_3,\ldots,\mu_5) \leq_{GDom} \partial^{( 6)}(\gamma/\lambda).$$
 We note that the above result will be still true if we apply row bumping algorithm on any row of $T$. That is  $ (\mu_2,\mu_3,\ldots,\mu_5) \leq_{GDom} \partial^{( i)}(\gamma/\lambda)$ for any $1\leq i\leq 7$. 
 
 \end{remark}

\vskip .1in
\subsection{ Proof of Theorem \ref{MainTheorem2}. } Now we provide the proof of  Theorem \ref{MainTheorem2} in which row bumping algorithm plays the crucial role.

\begin{proof}
	Given a skew partition $\gamma/\lambda$ with $n$ rows and $k$ columns, let $T$ be an LR-tableau of shape $\gamma/\lambda$ and content $\mu=(\mu_1,\mu_2,\ldots,\mu_i,\ldots )$   We already show in Proposition 4.1 that $ (\mu_{1},\mu_{2},\ldots) \leq_{GDom}    \partial^{(0)}(\gamma/\lambda)$. 
	
	Now for a sequence  $1\leq i_1<\ldots<i_j\leq n$,  we know that $ (\gamma/\lambda) ^{( i_1,\ldots,i_j)}$ be the skew  shape obtained by deleting the rows indexed by  $i_1,\ldots,i_j$ and concatenating remaining rows vertically. We want to show that  $$(\mu_{j+1},\mu_{j+2},\ldots,) \leq_{GDom}	\partial^{( i_1,\ldots,i_j)} (\gamma/\lambda) =^\leftarrow\uparrow (\gamma/\lambda^{( i_1,\ldots,i_j)}).$$
	
	First, assuming  that the size of the $ i_1$-th row  of $T$ is $r_1$,  observe that applying row bumping algorithm on the $ i_1$-th row  of $T$ results a semi-standard tableau $T^{(i_1)}$ of shape $(\gamma/\lambda)^{( i_1)}$ by removing $r_1$ labels $1$. So the content of $T^{(i_1)}$ is $(\mu_{1}-r_1,\mu_{2},\ldots)$ and the top cell of a column is labeled by $1$ or $2$. Observe that $T^{(i_1)}$ can be tiled with ballot-like  sequences $B_1,B_2,\ldots$ in which exactly  $r_1$ of  them  begins with label $2$. 
 
 Therefore we can still define row bumping algorithm on any row of $T^{(i_1)}$, since this algorithm is defined on the basis of ballot tiling of the given tableau. Since each label in this row lies in a different ballot-like sequence, the labels removed from $T^{(i_1)}$  are  $1$ or $2$, depending on whether  the ballot-like sequence of the label lying in the given row starts with $1$ or $2$. Hence in the resulting tableau, the top cell of a column is labeled by $1,2$ or $3$ and it can be tiled with ballot-like sequences starting with the labels $1,2$ or $3$.  
	
 In fact the above argument can be iterated for $j$ rows of $T$ indexed by $1\leq i_1<\ldots<i_j\leq n$ in the following manner: Observe that 
 $$T^{(i_1,i_2,\ldots,i_j)}:=(\ldots(T^{(i_1)})^{(i_2-1)}\ldots)^{(i_j-j+1)}$$
 is a semi-standard tableau that can be tiled with ballot-like sequence with the smallest label  $1,2,\ldots,j$ of $j+1$. Moreover its shape is $\gamma/\lambda^{(i_1,\ldots,i_j)}$ and  its content $(\mu_1',\ldots,\mu_j',\mu_{j+1},\mu_{j+2},\ldots)$ where $\mu_1+\ldots+\mu_j-\mu'_1-\ldots-\mu'_j$ is the sum of sizes of the removed rows indexed by $ i_1,\ldots,i_j$.

 Now let $T'$ be the semi-standard tableau  obtained by removing all the labels $1,2,\ldots,j$ from $T^{(i_1,i_2,\ldots,i_j)}$. 
 Now content of $T'$ is $(0,\ldots,0,\mu_{j+1},\mu_{j+2},\ldots)$ and it can be tiled with ballot-like sequences all of which starts with label $j+1$. Hence we can regard $T'$ as an LR-tableau with smallest label $j+1$. Denoting  its  shape  by $\gamma'/\lambda'$ we  then have by Proposition 4.1.,
 $$(\mu_{j+1},\mu_{j+2},\ldots)\leq_{GDom} \,^\leftarrow\uparrow(\gamma'/\lambda').$$ 
 On the other hand $ \gamma'/\lambda' \subset \gamma/\lambda^{(i_1,\ldots,i_j)}$  and hence  $$^\leftarrow\uparrow (\gamma'/\lambda') \leq_{GDom} \,^\leftarrow\uparrow 
 (\gamma/\lambda^{(i_1,\ldots,i_j)})=\partial^{( i_1,\ldots,i_j)} (\gamma/\lambda).$$
 Therefore we get $(\mu_{j+1},\mu_{j+2},\ldots)\leq_{GDom} \partial^{( i_1,\ldots,i_j)} (\gamma/\lambda)$ as required.

 \end{proof}

	\section{Final remarks}

It is well know fact that    $c_{\lambda,\mu}^{\gamma}=c_{\mu,\lambda}^{\gamma}= c_{\lambda^T,\mu^T}^{\gamma^T}=c_{\mu^T,\lambda^T}^{\gamma^T}$. Hence if   $c_{\lambda,\mu}^{\gamma}\not=0$,    we can conclude that for all $j\geq 0$ and for all $ 1\leq i_1<\ldots<i_j\leq n$, 
			\vskip.05in	
			\begin{itemize}
				\item[1.]  $\partial_{(j)}(\gamma/\lambda)\,\leq_{GDom}\,(\mu_{j+1},\mu_{j+2},\ldots) \leq_{GDom}   \partial^{(i_1,\ldots,i_j)}(\gamma/\lambda).$\vskip.05in	
				\item[2.]  
    $\partial_{(j)}(\gamma/\mu)\,\leq_{GDom}\,(\lambda_{j+1},\lambda_{j+2},\ldots) \leq_{GDom}   \partial^{(i_1,\ldots,i_j)}(\gamma/\mu).$\vskip.05in	
				\item[3.]  
    $\partial_{(j)}(\gamma^T/\lambda^T)\,\leq_{GDom}\,(\mu^T_{j+1},\mu^T_{j+2},\ldots) \leq_{GDom}   \partial^{(i_1,\ldots,i_j)}(\gamma^T/\lambda^T).$\vskip.05in	
				\item[4.]  
    $\partial_{(j)}(\gamma^T/\mu^T)\,\leq_{GDom}\,(\lambda^T_{j+1},\lambda^T_{j+2},\ldots) \leq_{GDom}   \partial^{(i_1,\ldots,i_j)}(\gamma^T/\mu^T).$\vskip.05in	
			\end{itemize} 

By using all upper and lower boundaries provided above, 	on some of the shapes $\gamma/\lambda$ which fits in $5\times 5$ and $6\times 6$ square we get the following tableau. Here first column represent size of $\gamma/\lambda$, second column represent actual number of triple $(\gamma,\mu,\lambda)$ with nonzero $c_{\lambda,\mu}^{\gamma}$. On the other hand last column represent the number triple $(\gamma,\mu,\lambda)$ which can  not be eliminated by our algorithm above although $c_{\lambda,\mu}^{\gamma}=0$. 

\begin{figure}[htp]
	\centering
	\includegraphics[width=12cm]{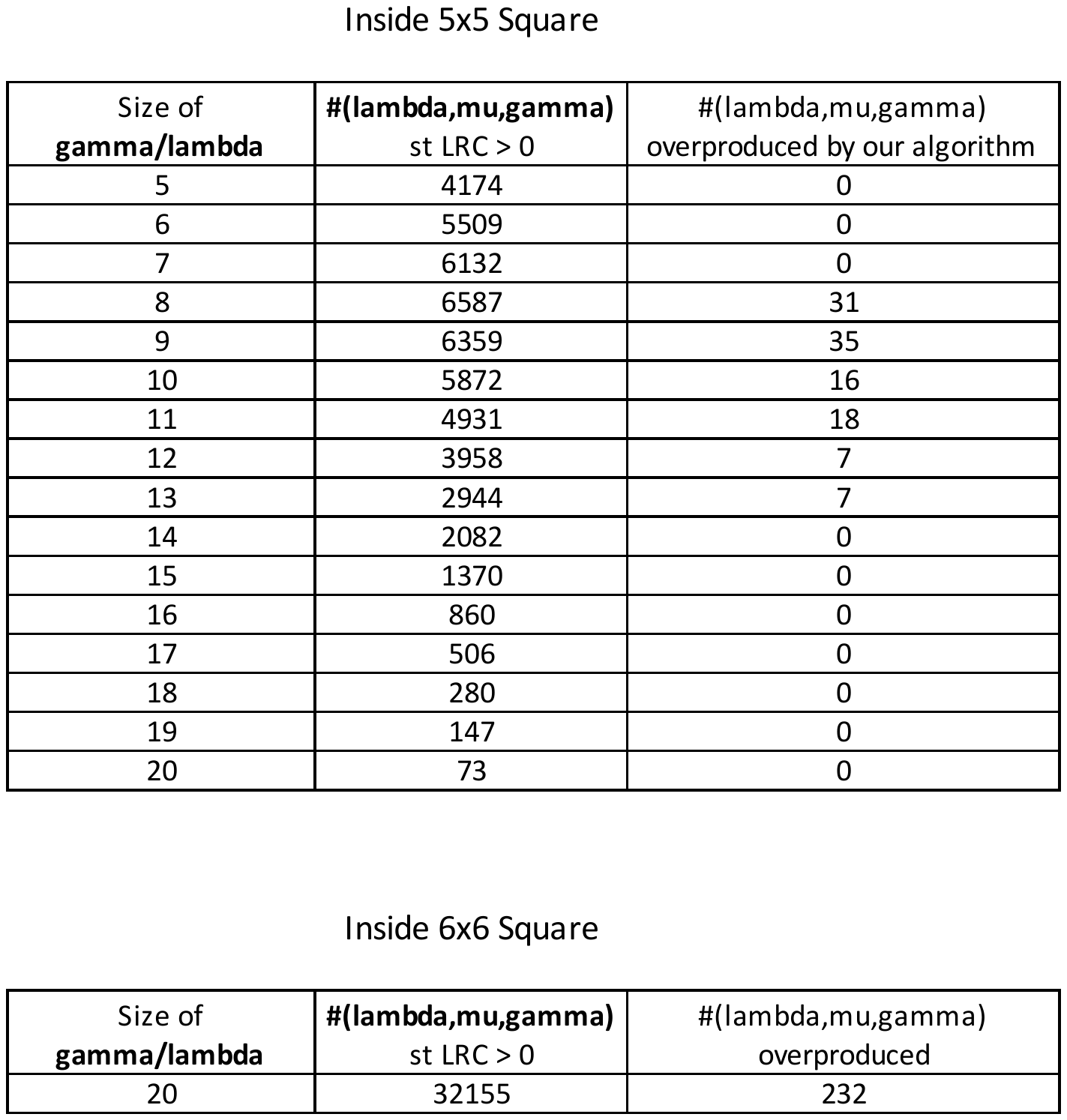}
\end{figure}

\end{document}